\documentclass[a4paper,12pt]{amsart}

\usepackage{amssymb,amsbsy,amsmath,amsfonts,amssymb,amscd}
\usepackage{latexsym}
\usepackage{graphics}
\usepackage{color}
\input xy
\xyoption{all}

\newcommand\sC{{\mathcal C}}

\newcommand\Ga{\Gamma}

\newcommand\ga{\gamma}

\newcommand{\CC}{\ensuremath{\mathbb{C}}}

\newcommand{\hol}{\ensuremath{\mathcal{O}}}

\newcommand{\HH}{\ensuremath{\mathbb{H}}}
\newcommand{\BB}{\ensuremath{\mathbb{B}}}
\newcommand{\PP}{\ensuremath{\mathbb{P}}}

\newcommand{\ra}{\ensuremath{\rightarrow}}

\def\eea{\end{eqnarray*}}
\def\bea{\begin{eqnarray*}}

\newcommand\dual{\mathrel{\raise3pt\hbox{$\underline{\mathrm{\thinspace d
\thinspace}}$}}}
\newcommand\qe{\ifhmode\unskip\nobreak\fi\quad $\Box$}       

\def\BOX{\hfill\lower.5\baselineskip\hbox{$\Box$}}

\newtheorem{theorem}{Theorem}[section]
\newtheorem{lemma}[theorem]{Lemma}
\newtheorem{proposition}[theorem]{Proposition}
\newtheorem{corollary}[theorem]{Corollary}

\newtheorem{claim}[theorem]{Claim}

\newtheorem{theo}{Theorem}[section]
\newtheorem{remarkk}[theo]{Remark}

\newtheorem{defin}[theo]{Definition}

\newtheorem{example}[theo]{Example}

\newcommand{\SSS}{\ensuremath{\mathcal{S}}}
\begin{document}

\title[Uniformization by a tube domain]{ A characterization of
varieties whose universal cover is the polydisk or a tube domain\footnote{AMS Classification: 32Q30, 32N05, 32M15, 32Q20, 32J25, 14C30, 14G35. }}
\author{ Fabrizio  Catanese - Antonio Jos\'e Di Scala}

\thanks{The present work took place in the realm of the DFG
Forschergruppe 790 "Classification of algebraic surfaces and compact
complex
manifolds". The visit of the second author to Bayreuth was supported
by the DFG FOR 790
The second author was also partially supported by GNSAGA (INdAM) and MIUR
(PRIN07, Differential Geometry and Global Analysis), Italy.}

\date{\today}

\maketitle

{\em  This article is  dedicated, with admiration \footnote{And with 
the friendship and gratitude of the first
author.},  to  Enrico Bombieri
on the occasion of his $70$-th birthday.}

\section{Introduction}

The uniformization theorem states that any complex manifold $C$ of
dimension 1 which is not  of special type (i.e., not $\PP^1$, $\CC$,
$\CC^*$,
or an elliptic curve) has as universal covering the unit disk $\BB_1
= \{ z  \in \CC | |z| < 1 \}$, which is biholomorphic to the upper
half plane
$\HH  = \{ z  \in \CC |Im (z) > 0 \}$.

A central problem in the theory of  complex manifolds has been the
one of determining the compact complex manifolds $X$ whose universal
covering
$\tilde{X}$ is biholomorphic to a bounded domain $\Omega\subset \CC^n$.

A first important restriction is given by theorems by   Siegel and
Kodaira,  extending to several variables  a  result of Poincar\'e,
and asserting
that necessarily such a manifold $X$ is projective and has ample
canonical divisor $K_X$ (see  \cite{Kodaira}, \cite{Kodaira-Morrow},
Theorem 8.4
page 144, where the Bergman metric is used, while the method of
Poincar\'e series is used in
\cite{Sie73},Theorem 3  page 117 , see also
\cite{Kollar}, Chapter 5).

In particular $X$ is a projective variety of general type  embedded
in projective space by a pluricanonical embedding associated to the
sections of
$\hol_X(mK_X)$ for large $m$.

This is a restriction on $X$, whereas a restriction on $\Omega$ is
given by  another  theorem of Siegel (\cite{Siegel}, cf. also
\cite{Kobayashi},
Theorem 6.2 \footnote{We are indebted to Pascal Dingoyan for 
providing this reference. }), asserting that $\Omega$ must be 
holomorphically
convex.

The problems which naturally come up  are then of two types:

{\bf Problem 1:} Given a
   bounded domain $\Omega \subset \CC^n$, when does there exist a
properly discontinuous group $\Ga \subset Aut (\Omega)$
   which acts freely on $\Omega$ and is cocompact  (i.e., is such  that
$ X = : \Omega / \Ga$ is a compact complex manifold with
   universal cover $ \cong \Omega$) ?

   The functions  on $\Omega$ which yield then a pluricanonical
embedding of $X$ are classically called automorphic functions,
   and in \cite[pag. 119]{Sie73} C.L. Siegel posed a second type of
problem writing:

\smallskip

\emph{... we have no method of deciding whether a given algebraic
variety of higher dimension can be uniformized by automorphic
functions.}

   A more specific  question than the one posed by Siegel is:

{\bf Problem 2:} Given a bounded domain $\Omega \subset \CC^n$, how
can we tell when a projective manifold $X$ with  ample canonical
divisor $K_X$
has $\Omega$ as  universal covering ?

Obviously an answer to the second problem presupposes an answer to
the first one.

For the first question it is natural to look at domains which have  a
big group of automorphisms, especially at {\bf bounded homogeneous
domains},
i.e., bounded domains such that the group $Aut(\Omega)$ of
biholomorphisms of $\Omega$ acts transitively.

And especially at the  {\bf bounded symmetric  domains}, the domains
such that for each point $p \in \Omega$ there is a symmetry at $p$ (an
automorphism $g$ with $ g(p)=p$ and $ (Dg)_p = - Identity$).

Bounded symmetric  domains were classified by Elie Cartan in
\cite{Cartan}, and  they are a finite number for each dimension $n$
(see also
\cite{Helgason2}, Theorem 7.1 page 383 and exercise D , pages 526-527
, and \cite{Roos} page 525 for a list of them).

Among them are the so called {\bf bounded symmetric domains of tube
type}, which are biholomorphic to a {\bf Tube domain}, a generalized
Siegel
upper halfspace
$ T_{\sC} =   \mathbb{V} \oplus  \sqrt -1  \sC$ where $\mathbb{V}$ is 
a real vector
space and $\sC \subset \mathbb{V}$ is a { \em symmetric cone}, i.e., a self dual homogeneous convex cone
containing
no full lines.

Borel proved in  \cite{Bo63} that for each bounded symmetric domain
$\Omega$ Problem 1 has a positive answer; and such a compact free
quotient $ X
= \Omega / \Ga$ is called a compact Clifford-Klein form  of the
symmetric domain $\Omega$.

Even if the bounded symmetric domains $\Omega$ are not the only ones
for which Problem 1 has a positive answer (i.e., such a compact
quotient $X$
exists),  Frankel proved in \cite{Frankel} that if $\Omega$ is a
bounded convex domain, and Problem 1 has a positive answer, then
$\Omega$ is a
bounded symmetric domain.

Another theorem of Frankel  (\cite{Fr95}){\footnote{We are indebted to Gang Tian for providing this reference} shows that $K_X$ ample implies the
splitting of  a finite unramified covering of $X$ as a product of a locally symmetric manifold and
a  {\em locally rigid} manifold, i.e., a manifold whose local group of isometries is discrete.

From this theorem of Frankel it follows that, if the universal cover is a bounded homogeneous domain,
then it must be a bounded symmetric domain. Since Theorem 1.1 of \cite{Fr95} asserts that 
 $ Aut ( \widetilde{X})$ acts as a group of isometries on $ \widetilde{X}$: therefore there can be no locally rigid factor if $ \widetilde{X}$
is a bounded homogeneous domain, hence $X$ is locally symmetric, so $ \widetilde{X}$ is a bounded symmetric domain.

Henceforth  we restrict our attention in this paper to Problem 2 for
the case where $\Omega$ is a bounded symmetric domain.

In this respect the first breakthrough, giving an answer to C.L.
Siegel's  question in an important special case, was based on the
theorems of
Aubin and Yau (see
\cite{Yau78}, \cite{Aubin}) showing the existence, on a projective
manifold with ample canonical divisor $K_X$, of a K\"ahler - Einstein
metric,
i.e. a  K\"ahler  metric
$\omega$ such that
$$ Ric (\omega) = - \omega . $$

This theorem is indeed the right substitute for the uniformization theorem in dimension  $n >1$.

Yau showed in fact  (\cite{Yau77}) that, for a projective manifold with
ample canonical divisor $K_X$,   the famous  Yau inequality is valid
$$ K_X^n \leq \frac{2(n + 1)}{n}  K_X^{n-2} c_2 (X),$$
  equality
holding if and only if the universal cover $ \widetilde{X}$ is the unit
ball $\BB_n$ in
$\CC^n$.

The uniformization theorems of Yau (\cite{Yau88}, \cite{Yau93}) for a
manifold $X$ with ample canonical bundle $K_X$ go in the direction of
providing further  answers to Siegel 's question, sketching  sufficient
(but not necessary) conditions in order that $\widetilde{X}$ be the 
product of  a
bounded symmetric domain with another manifold.

However Yau makes the unnecessary assumption that
   $\Omega^1_X$ splits as a direct sum $$\Omega^1_X= V_1 \oplus V_2
\oplus \cdots \oplus V_k,$$
    does not give an answer to the more specific Problem 2 and
moreover, as we shall show here, his conditions  for a summand $V_j$  apply only for
an irreducible factor of the universal cover which  is a ball or
a symmetric domain  of tube type.

A very readable exposition of Yau's  results, based on the concept of
stability of the cotangent bundle $\Omega^1_X$, is contained in the
first section of \cite{V-Z05}.

In the special case where  $\Omega^1_X$ splits as a sum of line
bundles it follows from Yau's theorem that $\widetilde{X}$ is the
polydisk
$\mathbb{H}^n$,
   where $n=\mathrm{dim}(X)$.

   The splitting of $\Omega^1_X$ as a sum of lines bundles is not a
necessary condition, even if  it does indeed
   hold on  a finite  unramified  covering $X' \rightarrow X$.

   The reason lies in the semidirect product (where $\SSS_n$ is the
symmetric group):

   $$ 1 \ra Aut(\HH)^n \ra Aut(\HH^n) \ra \SSS_n \ra 1.$$

    A necessary  condition for a compact complex manifold of dimension
$n$ to be uniformized by a polydisk was found in \cite{CaFr}, based
on the
consideration that the tensor (here $\odot$ denotes the symmetric product)
$$\tilde{\psi} =:  \frac{dz_1 \odot  \dots \odot dz_n }{ dz_1\wedge
\dots \wedge dz_n } $$

   is transformed by every automorphism $g$  into $ \sigma (g) \tilde{\psi}$,
   where $\sigma (g) = \pm 1$ is the signature of the permutation
corresponding to $g$.

Namely, the tensor $\tilde{\psi} $ descends to a so called semi
special tensor $\psi$ on $X$, which is simply a non zero section of
the sheaf
$S^n(\Omega^1_X)(-K_X) \otimes
\eta$, where $\eta$ is an invertible sheaf such that $\eta^2 \cong
\mathcal{O}_X$ (corresponding to the signature character).

   The necessary condition about the existence of a semi special tensor was
   proven, in dimension $n \leq 3$,  to be a sufficient condition for
$X$  to be uniformized by a polydisk ( \cite[Theorem 1.9.]{CaFr}).

Unfortunately, the above necessary  condition is not sufficient for
$n \geq 4$ (see \cite[Theorem 1.10.]{CaFr}).\\

Our first result in this paper is the following necessary and
sufficient condition for a compact complex manifold to be uniformized
by a polydisk.

\begin{theorem} \label{polydisc} Let $X$ be a compact complex
manifold of dimension $n$. Then the following two conditions:

\begin{itemize}

\item[(1)] $K_X$ is ample
\item[(2)] $X$ admits a semi special tensor $\psi \in
H^0(S^n(\Omega^1_X)(-K_X) \otimes \eta)$ such that, given any point
$p\in X$, the
corresponding hypersurface $F_p = :
\{\psi_p = 0 \}\subset \PP (TX_p)$ is  reduced

\end{itemize} hold if and only if $X \cong (\mathbb{H}^n)/ \Gamma$
(where $\Gamma$ is a cocompact discrete subgroup of
$\mathrm{Aut}(\mathbb{H}^n)$ acting freely ).
\end{theorem}

\begin{remarkk} The second condition is quite explicit, since it
amounts to verifying that the polynomial $\psi_p$, obtained evaluating
$\psi$ at the point $p$, is a square free polynomial: and to verify
this it suffices to use the G.C.D. of univariate polynomials.

\end{remarkk}

Our second and third results show that semispecial tensors, and a 
generalization of them, the
slope zero  tensors (see \cite{Bog} for the  related concepts of slope and stability) work out in a more
general setting, and give a \underline{necessary and sufficient
condition} for
a complex compact manifold $X$  to be uniformized by a bounded
symmetric domain of tube type.

Here a slope zero  tensor is a non zero section $\psi \in
H^0(S^{nm}(\Omega^1_X)(-m K_X) )$, where $m$ is a positive
integer.

\begin{theorem}\label{Main Theorem}

Let $X$ be a compact complex manifold of dimension $n$. Then the
following two conditions:

\begin{itemize}

\item[(1)] $K_X$ is ample
\item[(2)] $X$ admits a semi special tensor $\psi$;

\end{itemize}

hold if and only if $X \cong \Omega / \Gamma$ , where
$\Omega$ is a bounded symmetric domain of tube type with the special
property

(*) $\Omega$ is a product of irreducible bounded symmetric domains
$D_j$ of tube type whose rank $r_j$ divides the dimension $n_j$ of
$D_j$,

and $\Gamma$ is a cocompact discrete subgroup of
$\mathrm{Aut}(\Omega)$ acting freely.

Moreover, the degrees and the multiplicities of the irreducible
factors of the polynomial  $\psi_p$ determine uniquely the universal
covering
$\widetilde{X}=\Omega$.

\end{theorem}

\begin{theorem}\label{Main Theorem2}

Let $X$ be a compact complex manifold of dimension $n$. Then the
following two conditions:

\begin{itemize}

\item[(1)] $K_X$ is ample
\item[(2)] $X$ admits a slope zero   tensor $\psi \in
H^0(S^{mn}(\Omega^1_X)(-m K_X) )$, (here $m$ is a positive
integer);

\end{itemize} hold if and only if $X \cong \Omega / \Gamma$ , where
$\Omega$ is a bounded symmetric domain of tube type and $\Gamma$ is a
cocompact
discrete subgroup of
$\mathrm{Aut}(\Omega)$ acting freely.

Moreover, the degrees and the multiplicities of the irreducible
factors of the polynomial $\psi_p$ determine uniquely the universal
covering
$\widetilde{X}=\Omega$.

In particular, for $m=2$, we get that  the universal covering
$\widetilde{X}$ is a polydisk if and only if $\psi_p$ is the square of
a squarefree polynomial.

\end{theorem}

We obtain as  a corollary a simple proof of a variant  of  Kazhdan's
Theorem \cite{Kazh70} about the Galois conjugates of an arithmetic  projective
manifold
$X$. Namely, we have the following application.

\begin{corollary}\label{Kazhdan} Assume that $X$ is a projective manifold with
$K_X$  ample, and that
the universal covering $\tilde{X}$ is a bounded symmetric domain of tube type.

   Let $\sigma \in
\mathrm{Aut}(\mathbb{C})$ be an automorphism of $\mathbb{C}$.

Then the conjugate variety  $X^{\sigma}$ has universal covering
$\tilde{X^{\sigma}} \cong \tilde{X}$.
\end{corollary}

Our paper leaves two questions open:

\begin{enumerate}
\item
Is it possible (as in \cite{CaFr}) to remove the assumption that 
$K_X$ is ample, replacing
it by the condition that $X$ be of general type?

\item
Study  necessary and sufficient conditions for  the case where there 
are irreducible factors
which are bounded symmetric domains not of tube type: these should 
probably involve subbundles of higher rank
of the bundles $S^k(\Omega^1_X)(-mK_X)$.

\end{enumerate}

The paper is organized as follows: in   section 2 we recall a
result by Koranyi and Vagi which plays a central role for our
theorems, since it
determines the holonomy invariant hypersurfaces in the tangent space
to an irreducible symmetric bounded domain.

After this, in sections 3 and 4, we provide the proofs of our two
main theorems \ref{polydisc} and \ref{Main Theorem}, using the
existence of the
K\"ahler-Einstein metric, the classical theorems of De Rham and
Berger and the Bochner principle, in order to show the sufficiency of
the condition
of the existence of a semispecial tensor.

In section 4 we show that this condition is also necessary for every
bounded symmetric domain of tube type satisfying (*), thereby partly
generalizing  the result of Koranyi and Vagi (we prove invariance of
our tensor for the full group).

We conclude with the  Kazhdan type corollary \ref{Kazhdan}, a
couple of examples, and the proof of Theorem \ref{Main Theorem2}.

\section{Preliminaries}
\subsection{Symmetric bounded domains and its invariant polynomials}

Let $D \subset \mathbb{C}^n$ be a homogeneous bounded symmetric
domain in its circle realization around the origin $0 \in
\mathbb{C}^n$.

Let $K$ be the isotropy group of $D$ at the origin $0 \in
\mathbb{C}^n$, so that we have $ D =  G / K$.

   In \cite{KoVa79} polynomial $f \in \mathbb{C}[X_1,\dots,X_n]$ is said to be $K$-invariant
   (actually, it should be called {\em semi-invariant}) if there
    is a character $\chi \colon K \ra   \CC$ such that, for all $g \in K$,
$f(g X) =
\chi(g) f(X)$.

   Since $K$ is compact we have:  $|\chi(g)| = 1$.

Let $D = D_1 \times D_2$ be the decomposition of $D$ as a product of
two domains where $D_1$ is of tube type and $D_2$ has no irreducible
factor of
tube type. \\

\begin{theorem}\cite[Kor\'anyi-V\'agi]{KoVa79} \label{KV}

Let $D = D_1 \times D_2$ be the above decomposition and let moreover
   $$D_1 = D_{1,1} \times D_{1,2} \times \cdots \times D_{1,p}$$ be the
decomposition of $D_1$ as a product of irreducible tube type domains
$D_{1j},
\, \, \, (j=1,\cdots,p)$.

Then there exist, for each $j = 1, \dots p$, a unique
$K_j$-invariant polynomial $N_j (z_{1,j})$, where $K_j$ is the
isotropy subgroup of
$D_{1,j}$, such that:

   for all   $K$-invariant polynomial $f$   there exist a constant $c
\in \mathbb{C}$ and exponents $k_j$ with
   \begin{enumerate}
\item \[ f =  c \prod_{j=1}^p N_{j}^{k_j}, \, \, \, \] hence in particular
\item \[ f(z_1,z_2) = f(z_1) \, \, ,\] where $z_1$ denotes a vector
in  the domain $D_1$ and $z_2 \in D_2$.
\end{enumerate}
\end{theorem}

The above theorem follows from \cite{KoVa79} by taking into account
that a $K$-invariant polynomial is, up to a multiple, an \emph{inner
function},
i.e., a function such that $ | f(z) | = 1$ on the Shilov boundary of $D$.

This is so since the isotropy group $K$ acts transitively on the
Shilov boundary $S$ of $D$.

   It is very important to observe that the polynomials $N_j$ have 

degree equal to the $rank(D_j)$ of the irreducible domain $D_j$.

Here $rank(D_j)$ denotes  the dimension $r$ of the maximal totally 

geodesic embedded polydisc $\mathbb{H}^r \subset D_j$,
  or, 
equivalently, if $ D = G/K$, with $ G = Aut (D)^0$, $rank (D) = 
rank 
(G^{\CC}) =$ the dimension of the maximal
  algebraic torus contained 
in the complexification $G^{\CC}$.

   Therefore  $rank(D_j) \leq 
dim(D_j)$ and  equality holds if and 
only if $D_j = 
\mathbb{H}$.

The second part of the above theorem follows from part 
$(iii)$ in 
Theorem 3.3. of \cite[page 187]{KoVa79}.

The first part 
is contained in Lemma 2.5. and Lemma 2.3 of
\cite[pages 184,182]{KoVa79}.

The same result was rediscovered by Mok in \cite{Mok}.

\begin{remarkk}

The explicit form of the polynomial $N_j$ will be discussed in sections 4.2 and 4.3, 
where indeed we shall show that $N_j$ is (semi) invariant for the whole group $G$.

\end{remarkk}

\subsection{Irreducible symmetric domains of tube-type whose
dimension is divisible by its rank}

Recall the notation for the classical domains:

\begin{itemize}
\item
   $I_{n,p} $ is the domain $ D = \{ Z \in M_{n,p}(\mathbb{C}) :
\mathrm{I}_p - ^tZ \cdot \overline{Z} > 0 \}$.\\
\item
   $II_{n} $ is the intersection of the domain  $I_{n,n} $ with the
subspace of skew symmetric matrices.
   \item
   $III_{n} $ is instead the intersection of the domain  $I_{n,n} $
with the subspace of  symmetric matrices.
\end{itemize}

\begin{theorem}\label{tubes} Let $D$ be an irreducible symmetric
domain of tube-type. Let $d = dim(D)$ be the complex dimension of $D$
and $r$ its
rank.

If $d$ is multiple of $r$ then one of the following holds:\\

\begin{itemize}

\item[(i)] $D$ is of type $I_{n,n}$, $n \geq 1$. In this case $r=n$
and $d=n^2$,\\

\item[(ii)] $D$ is of type $II_{2k}$, $k \geq 1$. In this case $r=k$
and $d=k(2k-1)$,\\

\item[(iii)] $D$ is of type $III_{2k+1}$, $k \geq 0$. In this case
$r=2k+1$ and $d = (2k + 1)(k+1)$, \\

\item[(iv)] $D$ is of type $IV_{2k}$, $k \geq 2$. In this case $r=2$
and $d =2k$, \\

\item[(v)] $D$ is the exceptional domain of dimension $d=27$ and rank $r=3$.

\end{itemize}
\end{theorem}

\it Proof. \rm The proof follows from the classification of
irreducible bounded symmetric domains, see e.g. \cite[p. 525]{Roos}.

\section{Manifolds uniformized by a polydisk}

Here we prove Theorem \ref{polydisc}.

By a \emph{semi special tensor $\psi$ with reduced divisor}  we
mean, as in  \cite[Definition
1.3]{CaFr} ,
a semi special tensor $$\psi \in H^0 (S^n (\Omega^1_X)(- K_X) \otimes \eta)$$ such that
the homogeneous polynomial $\psi_p$, obtained evaluating the tensor on the fibre
over the point $ p \in X$
( $\psi_p$  is a polynomial of
degree $n$ on the tangent space $TX_p$),  is not divisible by a square.

Proposition 1.4. and its proof in \cite[page 6.]{CaFr} shows that, as
explained in the introduction,  (1) and (2) are necessary if $X \cong
(\mathbb{H}^n)/ \Gamma$ (where
$\Gamma$ is a cocompact discrete subgroup of
$\mathrm{Aut}(\mathbb{H}^n)$ acting freely ).\\

Assume now that (1) and (2) hold and let $\widetilde{X}$ be the
universal cover of $X$.

Proceeding as in \cite[page 160]{CaFr} the semispecial tensor $\psi$
pulls back to a special tensor $\tilde{\psi} = : \Psi $ on $
\tilde{X}$  which
is parallel with respect to the Levi-Civita connection associated to
the K\"ahler-Einstein metric.
   (this follows from the Bochner principle, see \cite{Koba80}, 
\cite{Yau88}, page 272 and \cite{Yau93}, page 479).\\

Fix a point $x \in \widetilde{X}$ and let $H_x \subset \mathrm{U}(T_x
\widetilde{X})$ be the restricted holonomy group with respect to  the
Levi-Civita connection associated  to the K\"ahler-Einstein metric.

   Since $\Psi$ is parallel there exists a degree $n$ polynomial $f :=
\psi_x$ on $T_x \widetilde{X}$
   such that \[ F_x = \{ v_x \in T_x \widetilde{X} :  \psi_x (v_x) = :
\Psi (x,v_x) = 0 \} \, \, \] is $H_x$-invariant.

This implies that $f = \psi_x$ is $H_x$-invariant in  the sense of
Koranyi-Vagi.

Notice that $f = \psi_x$ is not divisible by a square since $\psi_x$
has a reduced divisor $F_x$.\\

Since $\widetilde{X}$ has no flat De Rham factor (otherwise $X$ is
flat and  the canonical divisor $K_X$ cannot be ample) we use  the
second
author's   Proposition A.1 (appendix to \cite{CaFr}, page 178)
implying that there is a decomposition of the vector space $T_x
\widetilde{X}$ as
$T_x \widetilde{X} = V_1 \oplus V_2$ and where $f (v_1, v_2) =
f(v_1)$ depends only on the variable $v_1$.

 Moreover $V_1$ is the
tangent space
at the origin of a bounded symmetric domain $D \subset \mathbb{C}^m$
such that  the action of $H_x$ on $V_1$ is equal to the action of the
isotropy
group $K$ at the origin $0 \in \CC^m$ .\\

Let us use now Theorem \ref{KV} and notation therein. 

We obtain
  that $f$ splits as \[ f = c
\prod_{j=1}^p N_{j}^{\epsilon_j} \, \, \] where $\epsilon_j \in \{0,1
\}$. Then we get

\[ n = deg(f) = \sum_{j=1}^p \epsilon_j deg(N_{1j}) = \sum_{j=1}^p
\epsilon_j r_j \, \, . \]

We also have that

\[  \sum_{j=1}^p \epsilon_j r_j = n \geq m  = dim(D) = \sum_{j=1}^p
dim(D_{1j}) + dim(D_2) \geq \sum_{j=1}^p r_j + dim(D_2)\]

since $r_j \leq dim(D_{1j})$.

    We conclude that $\epsilon_j = 1$ $\forall j$, $n=m$, $p=n$,
$dim(D_2)= 0$, and
    moreover $dim(D_{1j}) = r_j = 1$ for $j=1,\cdots,n$.

     This shows that $H_x = K$ splits as $K = \mathrm{U}(1)^n$ and
completes the proof that $\widetilde{X}$ is a polydisk $\mathbb{H}^n$.

     \qed

\section{Manifolds uniformized by a tube domain}

Here we give the proof of Theorem \ref{Main Theorem}.

\subsection{Sufficient conditions}

We want here to show that if $K_X$ is ample, and $X$ admits a
semispecial tensor $\psi$, then the universal covering $\tilde{X}$
is a product of
irreducible symmetric domains of tube type whose rank divides the dimension.

We proceed as  for the proof of theorem \ref{polydisc}.

Namely, we write the universal cover  $\tilde{X}$, according to the
theorems of De Rham and Berger (see \cite{Berger} and also 
\cite{Olmos}),  as the product
$\tilde{X} = D_1 \times D_2 = D_1' \times D_1'' \times D_2 $ where
(since there are no flat factors, as already observed):

\begin{itemize}

\item
$D_2$ is the product of the irreducible factors of dimension $\geq 2$
for which the holonomy group is transitive (actually, it is the Unitary
group)
\item
$D_1$  is a bounded symmetric domain
\item
$D_1'$ is the product of all the irreducible bounded symmetric
domains of tube type.

\end{itemize}

Consider now the pull back tensor $\Psi = \tilde {\psi}$, and
consider coordinates $ (u,w,z) $ according to the product
decomposition $\tilde{X} =
D_1' \times D_1''
\times D_2 $.

Let $a = dim (D_1' ) , b =  dim (D_1'' ), r =  dim (D_2)$.

Then the tensor $\psi_x$ in a point $x$ can be written as
$$\psi_x =  f  (u,w,z) (du_1 \wedge \dots \wedge du_a)^{-1} \wedge
(dw_1 \wedge \dots \wedge dw_b)^{-1} \wedge  (dz_1 \wedge \dots \wedge
dz_r)^{-1}$$ and it is holonomy invariant.

By the same argument as in the previous section (Proposition A.1 of
the appendix to \cite{CaFr}, page 178, and the theorem of
Koranyi-Vagi) we have:
$$ f  (u,w,z) = f(u).$$

Write the (restricted) holonomy group as $K_1' \times K_1'' \times K_2$ and
observe that none of the subgroups  $K_2, K_1' , K_1''$ is 
contained 
in the
special unitary group, otherwise the K\"ahler-Einstein metric 
is 
Ricci flat in a certain direction, contradicting the ampleness of 
$K_X$.

Hence for instance $K_2$  acts non trivially on $ (dz_1 
\wedge \dots
\wedge dz_r)^{-1}$, while it acts trivially on $f$.

This would contradict the holonomy invariance of the tensor unless
there is no factor $D_2$. The same identical argument implies that
there is no
factor $D_1''$, hence $D$ is a product of  irreducible bounded
symmetric domains of tube type.

We write now accordingly $D$ as a product of such irreducible bounded
symmetric domains of tube type
$$  D = \prod_{j=1}^h \Omega_j $$ and we take variables $(z_1, \dots
, z_h)$ with $z_j \in \Omega_j$, and write, if

$z_j = (z_{j,1}, \dots , z_{j,n_j})$, and $n_j = dim (\Omega_j)$,

$$ dz_j^{top} = : dz_{j,1} \wedge  \dots \wedge dz_{j,n_j} .$$

By the theorem of Koranyi-Vagi, up to a constant we can write

$$\psi_x = N_1^{m_1}(z_1)  \dots N_h ^{m_h}(z_h)(dz_1^{top} \wedge
\dots    \wedge dz_h ^{top})^{-1}.$$

We impose invariance for each holonomy subgroup $K_j$.

We know that $K_j$ acts on $N_j(z_j)$ by a character $\chi_j(g)$, and
similarly  $K_j$ acts on $(dz_1^{top})$ by a character $\chi'_j$.

Recall that, $\Omega_j$ being a circular domain,  $K_j$ contains the
diagonal subgroup $S_j=  \{ e^{i \theta} I_{n_j}\}$.

Restricting to $S_j$ we see that, if  $\phi_j$  is the tautological
character, then
$\chi_j |_{S_j}  = \phi_j^{r_j} $ , $\chi' _j  |_{S_j}  =
\phi_j^{n_j} $, hence, by $S_j$ invariance,  we conclude that
$$ m_j r_j = n_j, \ \forall j=1, \dots h. $$

   We are done since we observe that the classification theorem
\ref{tubes} shows that the
   pair of integers $(r_j , n_j)$, under the condition $r_j  | n_j$,
completely determines
   the irreducible bounded symmetric domain of tube type $\Omega_j$.

\qed

\subsection{Necessary conditions}

As we observed in the introduction the ampleness of the canonical
line bundle $K_X$ is a result of Kodaira, i.e., condition $(1)$ is
necessary.\\

We shall give two proofs that condition $(2)$, i.e., the existence of
a semi special tensor,  is necessary.

Our first proof relies on the foundations of the theory of bounded
symmetric domains of tube type by means of their associated cones $\sC$ and
their Jordan
algebras, developed  for instance in
\cite{FaKo94}.

The second proof is a case by case computation which works  just for
the classical domains but provides an explicit expression for the semi special tensor.\\

Both proofs are based on the fact that, if $\Omega$ is a bounded 
symmetric domain,
and $$ \Omega = \Pi_{j=1}^h \Omega_j $$
is its decomposition as a product of irreducible bounded symmetric domains,
then we have a semidirect product
$$ 1 \ra  \Pi_{j=1}^h Aut (\Omega_j  ) \ra Aut(\Omega) \ra \SSS \ra 1$$
where $\SSS \subset \SSS_h$.

This follows from the fact that the De Rham decomposition of the universal cover 
of a complete Riemannian manifold is unique up to the ordering of the factors
(see \cite{KobNom}, Theorem 6.2 of Chapter IV).

As in the proof of Proposition 1.4 in \cite{CaFr} and by the above 
exact sequence it is enough to
construct, for each irreducible bounded symmetric domain of tube type
$D$, a
special tensor $\Psi$ invariant by the group of holomorphic
automorphisms $\mathrm{Aut}(D)$.

Then such a tensor $\Psi$ necessarily descends   to a semi special
tensor $\psi$ on any quotient $X$ of $\Omega$.

Let $D$ be an irreducible bounded symmetric domain of tube type.
Following \cite[Chapter X]{FaKo94} $D$
   is biholomorphic, via the Cayley map, to a tube domain $T_{\sC} =
\mathbb{V} + \mathrm{i} \sC$
    where $\mathbb{V}$ is a real finite dimensional vector space and
$\sC \subset \mathbb{V}$ is a so called symmetric cone.

  Both $D$ and
$T_{\sC}$ are open subsets of the Hermitian Jordan algebra
$\mathbb{V}_\mathbb{C} :=\mathbb{C} \otimes \mathbb{V}$ which is the
complexification of a simple Euclidean Jordan algebra whose real
vector space
is $\mathbb{V}$.

Let $\tilde{\psi}$ be the tensor defined as follows
\begin{equation}\label{psi} \tilde{\psi}:=
\frac{\mathrm{det}(\mathrm{d} z)^{\frac{n}{r}}}{\mathrm{K}}
\end{equation}

where $n = dim(D)$, $r$ is the rank of $D$,
$\mathrm{det}(\cdot)$ is defined in \cite[pag.29]{FaKo94} and
$\mathrm{K}$ is the complex
volume form of
$\mathbb{V}_\mathbb{C}$, i.e. a generator of
$\Lambda^n(\mathbb{V}_\mathbb{C})$.

  Notice that $\mathrm{det}(\cdot)$
is also denoted by
$\Delta(\cdot)$ and called the Koecher norm in \cite{KoVa79}. It is the same polynomial $N_j$ we encountered before.

Let $G(T_{\sC})$ be the group of biholomorphic maps of the tube
$T_{\sC}$.
\begin{lemma}\label{lemmaspecial}
$\tilde{\psi}$ is invariant by $G(T_{\sC})$.
\end{lemma}

\it Proof. \rm According to Theorem X.5.6 in \cite[p.207]{FaKo94} the
group $G(T_{\sC})$ is generated by the involution $j(z) :=
-z^{-1}$ and
the subgroups $G(\sC)$ and $N^+$. So it is enough to show that
$\tilde{\psi}$ is invariant by $j(z) := -z^{-1}$ and by the
subgroups
$G(\sC)$ and $N^+$.

That $\tilde{\psi}$ is invariant by the translations of $N^+$ is
obvious.

The invariance by $G(\sC)$ follows from Proposition
III.4.3 in
\cite[p.53]{FaKo94}.

To show that $j^* \tilde{\psi} =
\tilde{\psi}$ we will use the results in \cite[Chapter II]{FaKo94}
about the
   so called quadratic representation $P(\cdot)$, and also Lemma 1.1 
and Proposition 1.2
of \cite{ADO} stating the crucial properties:

\begin{itemize}
\item
$ P(x) (x^{-1} ) = x$
\item
$ P(x)^{-1}  = P(x^{-1}) $
\item
$ Dj (x) =  P(x)^{-1}$
\item
$ Det (P(x)) = det (x)^{\frac{2n}{r} }$
\item
$ det ( P(y) \cdot x) = (det y)^2 \cdot det \ x $.

\end{itemize}

   We have then:

$$ j^* \tilde{\psi} = \frac{\mathrm{det}(\mathrm{d}
j(z))^{\frac{n}{r}}}{j^*\mathrm{K}} =
\frac{\mathrm{det}( P(z)^{-1} \cdot 
\mathrm{d}z)^{\frac{n}{r}}}{j^*\mathrm{K}}= $$

$$ = \frac{\mathrm{det}( P(z^{-1}) \cdot 
\mathrm{d}z)^{\frac{n}{r}}}{j^*\mathrm{K}} =
\frac{( (\mathrm{det}z^{-1})^{2} \cdot 
\mathrm{det}(\mathrm{d}z))^{\frac{n}{r}}}{j^*\mathrm{K}}
=
\frac{(\mathrm{det}z)^{\frac{-2n}{r}} \cdot \mathrm{det}( 
\mathrm{d}z)^{\frac{n}{r}}}{j^*\mathrm{K}}=$$

$$=
\frac{
(\mathrm{det}z)^{\frac{-2n}{r}}.\mathrm{det}(\mathrm{d}z)^{\frac{n}{r}}}{\mathrm{Det}( 
P(z)^{-1})\mathrm{K}}  =
\frac{(\mathrm{det}z)^{\frac{-2n}{r}}.\mathrm{det}(\mathrm{d}z)^{\frac{n}{r}}}{(\mathrm{det} 
z)^{\frac{-2n}{r}}\mathrm{K}} = \tilde{\psi}.$$

  This
completes the proof of the claim.

\qed

\subsection{Necessary conditions found classically}

Here we construct explicitly $\tilde{\psi}$ for the classical bounded
symmetric domains of tube type.

We will follow the standard Elie
Cartan's notation.\\

{\bf Domains of type $I_{n,n}$}\\

The Cartan - Harish Chandra realization of $I_{n,n} := SU(n,n)/S(U(n) \times
U(n))$ is the domain $\Omega= \{ Z \in M_{n,n}(\mathbb{C}) :
\mathrm{I}_n - Z^t
\cdot \overline{Z} > 0
\}$.\\

To an element $\ga \in SU(n,n)$ corresponds the transformation
 $$\gamma(Z) = (AZ+B) \cdot (CZ + D)^{-1}.$$

As in
\cite[p. 174]{CaFr} the function $ \ga \mapsto \chi(\gamma) \in \mathbb{C}^*$
defined by the
equation:

\[ det(\mathrm{d} \gamma (Z)) = \chi(\gamma) \cdot det(CZ + D)^{-2}
\cdot det(\mathrm{d} Z) \]

is a character of $SU(n,n)$.

Indeed, if $\gamma' \in SU(n,n)$ is
another isometry, say 

$\gamma'(Z) = (A'Z+B') \cdot (C'Z + D')^{-1}$, then \[
det(\mathrm{d} (\gamma \cdot
\gamma')( Z )) = \chi(\gamma.\gamma') \cdot det((CA'+DC')Z +
CB'+DD')^{-2} \cdot det(\mathrm{d} Z) \] and by direct computation by
have

\[ \begin{aligned} det(\mathrm{d} (\gamma \cdot \gamma')( Z )) &=
\chi(\gamma) \cdot det(C\gamma'(Z) + D)^{-2} \cdot det(\mathrm{d}
\gamma'(z))\\
&= \chi(\gamma)
\chi(\gamma') \cdot det(C\gamma'(Z) + D)^{-2} \cdot det(C'Z +
D')^{-2} \cdot det(\mathrm{d} Z) \end{aligned}\]

so that  we only have to show that

\[ det((CA'+DC')Z + CB'+DD')^{-2} = det(C\gamma'(Z) + D)^{-2} \cdot
det(C'Z + D')^{-2} \]

which is equivalent to

\[ det((CA'+DC')Z + CB'+DD') = det(C\gamma'(Z) + D) \cdot det(C'Z + D') \]

but indeed

\[ \begin{aligned} det(C\gamma'(Z) + D) \cdot det(C'Z + D') & =
det(C(A'Z + B')\cdot(C'Z + D')^{-1} + D) \cdot det(C'Z + D') , \\
                   &=  det(C(A'Z + B')+ D(C'Z + D')) \\
                   &=  det((CA'+DC')Z + CB'+DD').\\
\end{aligned} .\]

This shows that $\chi(\gamma)$ is a character of $SU(n,n)$.Actually, any character of $SU(n,n)$ is trivial since
$SU(n,n)$ is a semisimple Lie group.

Hence the

\begin{claim}\label{Claim}
$\chi(\gamma) \equiv 1$, i.e. , the character $\chi$ is trivial.
\end{claim}  

Thus, for $\gamma \in SU(n,n) $,  we get the formula

\[ det(\mathrm{d} \gamma (Z)) =   det(CZ + D)^{-2} \cdot
det(\mathrm{d} Z) .\] \\

The Jacobian determinant of $\gamma$ is $det(CZ + D)^{-2n}$, i.e.
$\gamma^*\mathrm{K} = det(CZ + D)^{-2n} \mathrm{K}$, where
$\mathrm{K}$ is the
holomorphic volume form of
$I_{n,n}$.

Consider the tensor $\tilde{\psi}$ defined by \[ \tilde{\psi} =
\frac{det(\mathrm{d} Z)^n}{\mathrm{K}}\]

Then \[ \begin{aligned} \gamma^* \tilde{\psi} &=
\frac{det(\mathrm{d} \gamma Z)^n}{\gamma^* \mathrm{K}} \\ &=
\frac{(det(CZ + D)^{-2} \cdot
det(\mathrm{d} Z))^n}{det(CZ + D)^{-2n} \mathrm{K}} \\ &=
\frac{(det(CZ + D)^{-2})^n}{det(CZ + D)^{-2n}} \tilde{\psi} \\ &=
\tilde{\psi}
                                          \end{aligned}  \]

This shows that $\tilde{\psi}$ gives a special tensor which
descends to any Clifford-Klein form of $I_{n,n}$.\\
\vspace{1cm}

{\bf Domains of type $II_{2k}$.}\\

This is the subdomain of $I_{2k,2k}$ given by the skew-symmetric matrices.

Here $\tilde{\psi}$ is given by \[ \tilde{\psi} =
\frac{det(\mathrm{d} Z)^{\frac{2k-1}{2}}}{\mathrm{K}}\]

The Jacobian determinant of an isometry $\gamma$ is given by
\[ \gamma^* K = det(CZ + D)^{-(2k-1)} \mathrm{K} \, \, .\]

So \[ \begin{aligned} \gamma^* \tilde{\psi} &= \frac{det(\mathrm{d}
\gamma Z)^{\frac{2k-1}{2}}}{\gamma^* K} \\ &= \frac{(det(CZ + D)^{-2}
\cdot
det(\mathrm{d} Z))^{\frac{2k-1}{2}}}{det(CZ + D)^{-(2k-1)} K} \\ &=
\tilde{\psi}
                                          \end{aligned}  \]

This shows that $\tilde{\psi}$ gives a special tensor which
descends to any Clifford-Klein form of $II_{2k}$.\\

\vspace{1cm}

{\bf Domains of type $III_{2k+1}$.}\\

This is the subdomain of $I_{2k+1,2k+1}$ given by the symmetric matrices.

Here $\tilde{\psi}$ is given by \[ \tilde{\psi} =
\frac{det(\mathrm{d} Z)^{k+1}}{\mathrm{K}} \]

The Jacobian determinant of an isometry $\gamma$ is given by
\[ \gamma^* \mathrm{K}= det(CZ + D)^{-2(k+1)} \mathrm{K} \, \, .\]

So \[ \begin{aligned} \gamma^* \tilde{\psi} &= \frac{det(\mathrm{d}
\gamma Z)^{k+1}}{\gamma^* \mathrm{K}} \\ &= \frac{(det(CZ + D)^{-2}
\cdot
det(\mathrm{d} Z))^{k+1}}{det(CZ + D)^{-2(k+1)} \mathrm{K}} \\ &=
\tilde{\psi}
                                          \end{aligned}  \]

This shows that $\tilde{\psi}$ gives a special tensor which
descends to any Clifford-Klein form of $III_{2k+1}$.\\

\vspace{1cm}

{\bf Domains of type $IV_{2k}$, the so called Lie Balls.}\\

This domain admits a linear embedding into $I_{2^{2k},2^{2k}}$ via
Clifford algebras \cite[p.42]{Mc04}.\\

Here $\tilde{\psi}$ is given by \[ \tilde{\psi} =
\frac{det(\mathrm{d} Z)^{2k \cdot 2^{-2k}}}{\mathrm{K}}\]

The Jacobian determinant of an isometry $\gamma$ is given by
\[ \gamma^* \mathrm{K} = det(CZ + D)^{-2 \cdot 2k \cdot 2^{-2k}} \mathrm{K} \]

So \[ \begin{aligned} \gamma^* \tilde{\psi} &= \frac{det(\mathrm{d}
\gamma Z)^{2k \cdot 2^{-2k}}}{\gamma^* K} \\ &= \frac{(det(CZ + D)^{-2} \cdot
det(\mathrm{d} Z))^{2k \cdot 2^{-2k}}}{det(CZ + D)^{-2 \cdot 2k \cdot 2^{-2k}} K} \\ &=
\tilde{\psi}
                                          \end{aligned}  \]

This shows that $\tilde{\psi}$ gives a special tensor which
descends to any Clifford-Klein form of $IV_{2k}$.\\

\section{Proof of the Kazhdan's type corollary}

Consider the conjugate variety $X^{\sigma}$: since $K_X$ is ample we
may assume that $X$ is projectively embedded by $H^0 ( X, \hol_X (m
K_X)$.

$\sigma$ carries $X$ to
$X^{\sigma}$ and $K_X$ to $K_{ X^{\sigma}}$, hence also $X^{\sigma}$
has ample canonical divisor.

Consider a slope zero  tensor $\psi$ on $X$: then $\psi^{\sigma}$ is
also a slope zero   tensor, and moreover $\sigma$ sends the ring of
polynomial
functions on the tangent space
$TX_p$ to the corresponding ring of polynomial functions on the
tangent space $TX^{\sigma}_{\sigma(p)}$: hence the degrees and
multiplicities of
the irreducible factors of
$\psi_p$ are the same as   the degrees and multiplicities of the
irreducible factors of $\psi_{\sigma (p)}$.

We conclude then immediately by the last assertion of our main
theorems \ref{Main Theorem} and \ref{Main Theorem2} that  the 
universal covering of
$X^{\sigma}$ is $\tilde{X}$.

\qed

\section{Examples}

Assume that the polynomial $\psi_p$ associated to a semi special 
tensor is a square $\psi_p  = N^2$, where $N$
is irreducible (the more general case where $N$ is square free
follows then
right away).

Then the universal covering $\tilde{X}$ is an irreducible symmetric
tube domain such that $d/r = 2$.

It follows from  Theorem \ref{tubes} that $\tilde{X}$ is either
$I_{2,2}$ or $III_3$. In particular  $X$ has either dimension 4 or 6.

\begin{proposition}\label{square} 
Assume that $K_X$ is ample and $X$ admits a
semispecial tensor $\psi$.

If the multiplicities of the divisor associated to $f= : \psi_p$ are
at most $2$ then $\tilde{X}$ is a product of 1-dimensional disks, of
domains
of type $I_{2,2}$ or of type $III_3$.

Moreover, if all multiplicities are $2$ then the number of factors of
$f$ and  the dimension $n$ of $X$ determine $\tilde{X}$.
\end{proposition}

\it Proof. \rm The hypotheses imply that the polynomial  $f =
\psi_p$ can be factorized as
\[ f = c \prod_{j=1}^p N_{j}^{e_j} \, \, \]

where $e_j \leq 2$. If $e_j = 1$ then the corresponding factor is a disk.

If $e_j = 2$ by the previous observation the corresponding factor is
is either $I_{2,2}$ or $III_3$, and this shows the first assertion.

The hypothesis of the second statement is that
\[ f = c \prod_{j=1}^p N_{j}^{2} \, \, \]

Let us denote by  $a$ the number of times that $I_{2,2}$ occurs in
$\tilde{X}$ and by $b$
   the number of times that $III_{3}$ occurs in $\tilde{X}$.

   Then
\[ \begin{cases} 4a + 6b = n = dim(X) \\ a + b = p
\end{cases} \]

Hence, knowing   $p$ and $n$,  we know $a,b$ and also  $\tilde{X}$.

\qed

\section{Slope zero tensors of higher degree}

Let's treat first the case where  $\tilde{X}$ is an irreducible symmetric bounded domain of tube
type of dimension $n$ and rank $r$, but where we consider more generally  the sheaf
$S^k(\Omega^1_{\tilde{X}})(-m K_{\tilde{X}})$,  $k,m$ being  positive
integers.

Assume that there exists a  tensor $\tilde{\psi} \in
H^0(\tilde{X},S^k(\Omega^1_{\tilde{X}})(-m K_{\tilde{X}}))$ invariant
by the full automorphism group $Aut(\tilde{X})$.

  Then by the theorem
of
Koranyi-Vagi
\begin{equation} \tilde{\psi}_x = N^{a}(z)( dz^{top})^{-m}
\end{equation}
   and
\begin{equation}\label{mnra} k = m \cdot n = r \cdot a
\end{equation} since $\tilde{\psi}$ is invariant by the diagonal
subgroup $S = \{ e^{i \theta} I_{n}\}$.

  Conversely, if condition
(\ref{mnra})
holds then $\tilde{\psi}$ is invariant by the full group of
automorphisms (the proof is the same as in \ref{lemmaspecial}), hence
$\tilde{\psi}$
descends to any Clifford-Klein form X of $\tilde{X}$: providing a
section $\psi$ of the sheaf $S^k(\Omega^1_{X})(-m K_{X})$.

Let now $\tilde{X}$ be the product $\Omega_1 \times \cdots \times
\Omega_h$ of the irreducible symmetric bounded domains of tube type of
dimension
$n_j$ and rank $r_j$,
$j=1,\cdots,h$.

  If $\tilde{\psi} \in
H^0(\tilde{X},S^k(\Omega^1_{\tilde{X}})(-m K_{\tilde{X}}))$ is
invariant by $\mathrm{Aut}(\Omega_1) \times
\cdots \times
\mathrm{Aut}(\Omega_h)$ then
\begin{equation}
\psi_x = N_1^{a_1}(z_1)  \dots N_h ^{a_h}(z_h)( dz_1^{top} \wedge
\dots    \wedge dz_h ^{top})^{-m}.
\end{equation}  $k = m \cdot n$ and $a_j \cdot r_j = m \cdot n_j$ 
for $j=1,\cdots,h$.

Conversely, if the numerical conditions $a_j \cdot r_j = m \cdot n_j$ hold for
$j=1,\cdots,h$, then the above formula for $\psi_x$ defines a
section of the
sheaf
$S^k(\Omega^1_{\tilde{X}})(-m K_{\tilde{X}})$ invariant by
$\mathrm{Aut}(\Omega_1) \times \cdots \times
\mathrm{Aut}(\Omega_h)$.\\

Now notice that,  for any product $\Omega_1 \times \cdots \times
\Omega_h$ of  irreducible symmetric bounded domains  of tube type of
dimension
$n_j$ and rank $r_j$,  we can {\bf always} find integers
$m,a_1,\cdots,a_j$ such that the numerical conditions $a_j \cdot r_j =
m \cdot n_j$ hold for
$j=1,\cdots,h$.

By using the 2-torsion invertible sheaf $\eta$ corresponding to the signature (of $\SSS 
\subset \SSS_h$) we get a non zero section 
$$ \psi \in H^0 (S^{mn}(\Omega^1_{\tilde{X}})(-m K_{\tilde{X}}) \otimes \eta).$$

If $\eta$ is nontrivial, replace $\psi$ by $\psi^2$: we  obtain in this way  a slope zero tensor.

Hence

{ \bf Theorem \ref{Main Theorem2}}

{\em Let $X$ be a compact complex manifold of dimension $n$.

Then the
following two conditions:

\begin{itemize}

\item[(1)] $K_X$ is ample
\item[(2)] $X$ admits a slope zero  tensor $\psi \in
H^0(S^{mn}(\Omega^1_X)(-mK_X) )$, where $m$ is a  positive
integer;

\end{itemize} hold if and only if $X \cong \Omega / \Gamma$ , where
$\Omega$ is a bounded symmetric domain of tube type and $\Gamma$ is a
cocompact
discrete subgroup of
$\mathrm{Aut}(\Omega)$ acting freely.

Moreover, the degrees and the multiplicities of the irreducible
factors of the polynomial $\psi_p$ determine uniquely the universal
covering
$\widetilde{X}=\Omega$.

In particular, for $m=2$, we get that  the universal covering
$\widetilde{X}$ is a polydisk if and only if $\psi_p$ is the square of
a squarefree polynomial.}

\bigskip

The proof is identical to the proof of Theorem \ref{Main Theorem}
taken into account the observation made above (for the existence 
part) that it is
possible to
find the numbers
$m,a_1,\cdots,a_j$ such that the numerical conditions $a_j \cdot r_j =
m \cdot n_j$ holds for $j=1,\cdots,h$.\\

Here is one more example.

Let $X$ be a compact $3$-dimensional complex
manifold with $K_X$ ample and such that $\psi \in
H^0(S^{3m}(\Omega^1_X)(-mK_X) \otimes
\eta)$.

Then either
$\tilde{X}=\mathbb{H} \times \mathbb{H} \times \mathbb{H}$ or
$\tilde{X}$ is the Lie ball, i.e., the domain of type $IV$ and
dimension $3$.

  In this
last case the sheaf
$S^6(\Omega^1_X)(-2K_X)$ has a section.

Notice that the rank=2 does not
divide the dimension=3 and that the divisor of the section is not
reduced.

{\bf Acknowledgement:}

We would like to thank Marco Franciosi for interesting conversations
which led to our present cooperation.

\medskip
\noindent {\bf Authors '  Addresses:}\\
\noindent Fabrizio Catanese, \\ Lehrstuhl Mathematik VIII, Mathematisches
Institut der \\Universit\"at Bayreuth\\ NW II,  Universit\"atsstr. 30,
95447
Bayreuth, Germany.\\
            email:          fabrizio.catanese@uni-bayreuth.de\\

\noindent Antonio Jos\'e Di Scala, \\ Dipartimento di Matematica,
Politecnico di Torino, \\ Corso Duca degli Abruzzi 24, 10129 Torino,
Italy. \\
    email:     antonio.discala@polito.it \\

\end{document}